\documentclass[twocolumn]{autart_edited}    

\usepackage{graphicx}          

\usepackage{amsmath,amssymb,amsfonts}
\usepackage{graphicx}
\usepackage{textcomp}
\usepackage{algpseudocode}
\usepackage{bm}
\usepackage{color}

\newtheorem{theorem}{Theorem}
\newtheorem{corollary}{Corollary}
\newtheorem{proposition}{Proposition}
\newtheorem{definition}{Definition}
\newtheorem{remark}{Remark}
\newtheorem{LP}{Linear Programming Problem}

\newcommand{\co}{\mathsf{co}}
\newcommand{\ee}{\epsilon}

\begin{document}

\begin{frontmatter}

\title{Control Design with Guaranteed Transient Performance: an Approach with Polyhedral Target Tubes\thanksref{footnoteinfo}} 
\thanks[footnoteinfo]{© 2020. This manuscript version is made available under the CC-BY-NC-ND 4.0 license http://creativecommons.org/licenses/by-nc-nd/4.0/. DOI: 10.1016/j.automatica.2020.109097. Corresponding author: Willem Esterhuizen.}

\author[WE_affiliation]{Willem Esterhuizen}\ead{willem.esterhuizen@etit.tu-chemnitz.de},    
\author[QW_affiliation]{Qing-Guo Wang}\ead{wangq@uj.ac.za},               

\address[WE_affiliation]{Technische Universit\"{a}t Chemnitz, Automatic Control \& System Dynamics Laboratory, Germany and Institute for Intelligent Systems, Faculty of Engineering and the Built Environment, University of Johannesburg, South Africa}  
\address[QW_affiliation]{Institute for Intelligent Systems, Faculty of Engineering and the Built Environment, University of Johannesburg, South Africa}             
          

\begin{abstract}                          
In this paper a novel approach is presented for control design with guaranteed transient performance for multiple-input multiple-output discrete-time linear polytopic difference inclusions. We establish a theorem that gives necessary and sufficient conditions for the state to evolve from one polyhedral subset of the state-space to another. Then we present an algorithm which constructs a time-varying output feedback law which guarantees that the state evolves within a time-varying polyhedral \emph{target-tube} specifying the system's desired transient performance. We present generalisations involving constraints on the control, and a bounded additive disturbance term. Our formulation is very general and includes reference tracking with any desired transient behaviour in the face of disturbances, as specified, for example, by the most popular step response specifications. The approach is demonstrated by an example involving the control of water levels in two coupled tanks.
\end{abstract}

\end{frontmatter}

\section{Introduction}

Research in control theory is often concerned with stability without regarding the transient behaviour of a system. However, in industry this aspect is very important: a stable system with big transient errors or a very long settling time is obviously undesirable and very likely to be put aside by engineers. PID has enjoyed much success in industry not only because it is relatively simple, but also because the design often addresses performance specifications, such as rise-time, settling-time, over-shoot and steady-state error, see \cite{Wang97}. However, one of PID's shortcomings is the fact that it does not \emph{explicitly} take constraints and performance requirements into account, resulting in engineers often resorting to their experience and trial-and-error.

The problem of designing control systems capable of shaping transient performance has received the attention of a number of researchers, and we present a summary of the problem's history. Typical adaptive control schemes are plagued by unacceptable transients and it is in this field that, to our knowledge, one of the first papers on shaping a system's transient response according to explicit performance specifications, see \cite{Miller91}, was produced. Concentrating on single-input single-output (SISO) systems that are minimum phase, the authors approach the problem by dynamically adjusting the controller's feedback gains. Other works in adaptive control that are concerned with transient performance usually present results in the form of guaranteed bounds on the evolution of the state, output, or control signals, see for example \cite{Zang94,KRSTIC1993}. Note, however, that these works do not consider the derivation of control laws such that explicit performance requirements are satisfied.

Funnel control, introduced in \cite{Ilchmann:2002}, approaches the problem by specifying time-varying constraints (the ``funnel'') on the output and letting the control magnitude be dependent on the distance of the output to the funnel boundary. The theory is applicable to a great diversity of dynamical systems of known relative degree that satisfy a ``high frequency gain condition'', and which have the same number of inputs as outputs. Later works, see \cite{Ilchmann:2010:SISO,Ilchmann:2010:MIMO}, consider constrained inputs in the formulation; and \cite{Liberzon:2013} considers a bang-bang implementation. Funnel control has found application in many fields, see for example \cite{POMPRAPA_2015,Hackl_2011,Hackl_2013_b}. The reader may also refer to the references \cite{Hackl_2017,Ilchmann_Mueller:2009,Ilchmann:2007_2}.

Another approach to shaping transients is presented in \cite{BechlioulisR08MIMO} and \cite{Bechlioulis:2009SISO} for feedback linearisable nonlinear systems, and systems in strict-feedback, respectively. The authors specify performance requirements through time-varying constraints on the state and introduce a transformation that recasts this problem into an unconstrained one, the stability of which leads to a solution of the original performance problem. Though the method is applicable to unknown nonlinear systems with mild assumptions, one still needs to solve a new stabilisation problem. More recently, the paper \cite{Bechlioulis2014} introduces a similar, but simpler, idea to systems in pure feedback form.


Then there are approaches that may be described as \emph{target-tube} problems (introduced in \cite{BERTSEKAS_1971} and \cite{Schweppe_71}), where one specifies a time-varying set, called the target-tube, over a finite horizon and then seeks a control law that keeps the state in this tube. In \cite{BERTSEKAS_1971} the authors showed that under the assumption that the target-tube and control constraints are ellipsoidal, a linear time-varying feedback can be constructed backwards in time for discrete-time linear systems with an additive disturbance term. They also presented results for polyhedral target tubes, under the assumption of an invertible system matrix. The work of \cite{Kurzhanski_1997,KURZHANSKIY_2011} and \cite[Ch.8]{Kurzhanski_2014}, addresses the problem using an ``ellipsoidal calculus''. The approach involves finding a problem's ``solvability set'' through the solution of a ``funnel equation'' and at every time instant selecting a control value that results in the state evolving towards the solvability set's interior. Under the assumptions of linear systems and ellipsoidal constraint sets the solvability set may be obtained from the solutions of differential equations. The papers \cite{Kerrigan_2004} and \cite{KERRIGAN_2006} consider linear discrete-time systems subjected to mixed polyhedral constraints and a bounded additive disturbance term, with the goal of finding a sequence of \textit{disturbance feedback laws}. They show that the set of feedbacks is convex, but the approach requires exact knowledge of all past states and inputs and an accurate model of the system. The work \cite{Kostousova_2018} considers the problem for a class of uncertain discrete-time systems under the assumption of parallelepipeds in the problem data. An iterative algorithm to solve the control synthesis problem is then presented.

In this paper we present a novel target-tube approach to ensuring desired transient performance. We build on our previous conference paper, see \cite{ESTERHUIZEN_2017}, where some of the ideas in this work were introduced, and present a much improved and more complete approach, along with some results on the existence of a solution to our introduced algorithm. The contributions are summarised as follows.

\begin{enumerate}
	\item We consider discrete-time \emph{linear polytopic difference inclusions}. This type of model, which can be interpreted as an uncertain linear model, also often appears in the literature as a surrogate model for some nonlinear systems. Examples include \emph{linear parameter varying systems}, see \cite[Ch. 7]{Blanchini:2015} and \cite{Toth_2010}, and fuzzy nonlinear system modelling, \cite{Tanaka_2001}. Most other works that address the target-tube problem assume an exact linear time-invariant model, possibly with a bounded additive disturbance term.
	\item We require the output to be contained in a time-varying \emph{polyhedral} target-tube, which specifies the desired performance characteristics. Polyhedral constraints are often more natural to impose on physical systems in comparison with, for example, time-varying ellipsoidal constraint sets, which the vast majority of other works on target-tubes consider.
	\item We present an algorithm where we start with the target set, $\mathcal{X}_T$, at the end of the horizon and find a sequence of output feedbacks, of the form $u(k) = F(k)y(k)$ with $F(k)$ a matrix for every $k$, \textit{backwards} in time in a dynamic programming fashion, as in \cite{BERTSEKAS_1971}, along with an initial set, $\mathcal{X}_0$, such that for \textit{every} initial condition $x_0\in\mathcal{X}_0$ the feedback results in the state evolving within the target-tube until it reaches the target set. This is an improvement over the method used in \cite{ESTERHUIZEN_2017} where the feedback was found forwards in time, which made the constructed sets susceptible to the wrapping effect. To our knowledge, the current paper, as well as the paper \cite{ESTERHUIZEN_2017}, are the first that attempt to find a \emph{time-varying linear state/output feedback} under the assumption of a polyhedral target-tube.
	\item We emphasise that in our approach the sequence of feedbacks is found \emph{off-line} via the solution of linear programming problems, and that there are \emph{no online optimisation problems} that need to be solved. Thus, our approach shares no similarity with model predictive control.
	\item We generalise the results so that one may include a \textit{bounded additive disturbance term} in the dynamics, as well as \textit{polytopic constraints on the control}.
\end{enumerate}

To arrive at our results we exploit ideas that appear in \emph{set-theoretic methods}. In that context, when studying the performance of systems, one is usually concerned with concepts such as decay-rate, transient estimates and reachable sets, see for example \cite{Blanchini:2015}, \cite{Boyd:94} and \cite{Girard:2006}, and the references therein. We adapt a well-known result that appears in \cite{Bitsoris:Jun1988} and \cite{Blanchini:2015}, which is related to polyhedral invariant sets, to produce a theorem with conditions that guarantee that the state evolves from one polyhedral set to another. We then use this theorem to arrive at an algorithm that constructs a feedback backwards in time, as mentioned in point (3) above.

The paper is organised as follows: in Section~\ref{sec:problem_formulation} we specify the system under investigation and the problem we wish to solve. In Section~\ref{sec:reach_sets} we present an overview of reachable sets and state our main tool: a theorem that gives necessary and sufficient conditions for the state to evolve from one polyhedral set to another. Section~\ref{sec:algorithm} presents an algorithm that uses Theorem~1 to construct a time-varying feedback that solves our problem. In Section~\ref{sec:generalisations} we generalise the results to systems with a bounded additive disturbance term and polytopic constraints on the control. Section~\ref{sec:stab} is dedicated to a discussion on stability issues. In Section~\ref{sec:example} we present an example, involving the control of fluid levels in two coupled tanks via output feedback. Finally, we conclude with Section~\ref{sec:conclusion}.

\subsection*{\textbf{Notation}}\label{sec:notation}

If $M$ is a matrix, $M^{[i]}$ refers to its $i$-th row. The notation $M\geq 0$ means that every element of the matrix $M$ is nonnegative. ($M\geq 0$ does \emph{not} mean $M$ is positive-semidefinite.) If $r$ is an $n$-dimensional vector, then $r_i$ refers to its $i$-th coordinate, and $r\geq 0$, $r\leq 0$ and $r\neq 0$ indicates that $r_i\geq 0$, $r_i\leq 0$ and $r_i\neq 0$, respectively, for $i\in\{1,\dots,n\}$. If both $r$ and $s$ are $n$-dimensional vectors, then the notation $r\leq s$ is to be interpreted element-wise. The notation $r^T$ indicates the transpose of the vector $r$. A column vector of appropriate dimension with all its elements equal to one is given by \textbf{1}. A matrix, of appropriate dimensions, with all its elements equal to 0 is given by $\mathbf{0}$. $I_{n\times n}$ refers to the $n$ dimensional identity matrix. The notation $\mathbb{R}^n$ refers to $n$-dimensional Euclidean space; $\mathbb{R}^{m\times n}$ to the set of all matrices with $m$ rows and $n$ columns, and with real entries; and $\mathbb{Z}_{\geq 0}$ to the set of nonnegative integers. The acronym \emph{s.t.} stands for ``subject to''. A \emph{polyhedral set} is specified by $\mathcal{P}(M,m) = \{x\in\mathbb{R}^n : Mx \leq m\}$ where $M\in\mathbb{R}^{p\times q}$, with \emph{not} all the elements of the row vector $M^{[i]}$ equal to 0; and $m\in\mathbb{R}^p$. 
If $\mathcal{S}\subset \mathbb{R}^n$,  $TS\triangleq \{Ts:s\in S\}$ and  $\beta\mathcal{S} = \{\beta s, s\in \mathcal{S}\}$, where $T$ is a linear transformation and $\beta\geq 0$ a real number. Given two subsets of $\mathbb{R}^n$, $S_1$ and $S_2$, the \emph{Minkowski sum} is given by $S_1\oplus S_2 \triangleq\{s_1 + s_2: s_1\in S_1, s_2 \in S_2\}$. A \emph{C-set}, \cite{Blanchini:2015}, is a convex and compact subset of $\mathbb{R}^n$ containing the origin as an interior point. The C-set $\mathcal{S}$ is said to be \emph{$\eta$-contractive}, \cite{Blanchini:2015}, for the system $x(k+1) = Ax(k)$ provided that for every $x\in\mathcal{S}$ we have $\mathcal{M}_{\mathcal{S}}(Ax) \leq \eta$, $0\leq\eta < 1$. The notation $\co\{M_i\} = \{M = \sum_{i=1}^s \beta_i M_i:\beta_i\geq 0, \sum_{i=1}^s\beta_i = 1\}$, where the summation is done element-wise.

\section{Problem Formulation}\label{sec:problem_formulation}
We consider linear polytopic difference inclusions:
\begin{align}
x(k + 1)  & =   \bm{A}(k)x(k) + \bm{B}(k)u(k)\label{eq:systemeq1}\\
y(k) & = Cx(k), \label{eq:systemeq2}
\end{align}
where $k = 0,1,\dots,K$ is the time index and $K\in \mathbb{Z}_{\geq 0}$ specifies the time horizon; $x(k)\in \mathbb{R}^n$ is the state, $y(k)\in\mathbb{R}^r$ is the output, and $u(k)\in\mathbb{R}^m$ is the control. The concatenated matrix $[\bm{A}(k)\,\, \bm{B}(k)]$ is assumed to be ``polytopic'', that is, $[\bm{A}(k)\,\,\bm{B}(k)] \in \co\{[A_i\,\,B_i]\}_{i=1}^s$ for every $k$, and $C$ is assumed to be constant.

We now introduce two polyhedral sets: with $Q_T\in\mathbb{R}^{q_T \times n}$ a constant matrix,  $\mathcal{X}_T\subset \mathbb{R}^n$ is the \emph{target} set:
\[
\mathcal{X}_T = \mathcal{P}(Q_T,\psi_T),
\]
with $\psi_T$ having compatible dimension. With $Q(k)\in\mathbb{R}^{q(k)\times n}$ (i.e., the number of rows of the matrix $Q(k)$ may vary in time), and $\phi(k)$ a vector of compatible dimension for every $k$, the performance requirements of the system are specified with a time-varying polyhedral set, called the \textit{target-tube}:
\[
\mathcal{H}(k) = \mathcal{P}(Q(k),\phi(k)),\quad k \in\{0,\dots,K\},
\]
with $\mathcal{H}(K) = \mathcal{X}_T$. We make the following assumption, which is introduced in connection with the existence of solutions to our algorithm, which is to follow:
\begin{description}
	\item[(A1)] $\mathcal{H}(k)$ is bounded for every $k\in\{0,1,\dots,K\}$.
\end{description}

\begin{figure}[thpb]
	\begin{center}
		\includegraphics[width=0.8\columnwidth]{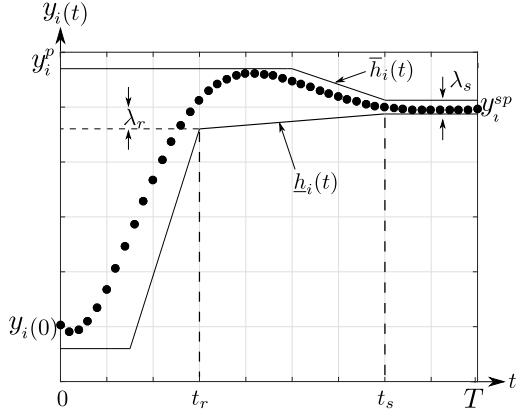}
		\caption{Transient performance specifications for a step set-point change. The set-point is specified by $y_i^{sp}$, the steady-state error by $\lambda_s$, the $\lambda_r$-rise-time by $t_r$ and the peak over-shoot by $y_i^p$.}
		\label{fig:perf_func}
	\end{center}
\end{figure}

We point out that our formulation allows for the specification of transient performance for the state \emph{or} output, and that one may opt to only specify the performance of \emph{some} of the coordinates. As an example, suppose that for a planar system there are two outputs, $y_1 = x_1 + x_2$ and $y_2 = x_2 - x_1$, subjected to constant performance constraints: $|y_1|\leq 1$ and $|y_2|\leq 1$ for all $k$. We would then have $Q(k) \equiv \left(
\begin{array}{cccc}
1 & -1 & -1 & 1\\
1 & -1 & 1 & -1
\end{array}
\right)^T\,\, \phi(k) \equiv \left(\begin{array}{cccc}
1 &1 &1 & 1
\end{array}
\right)^T.
$ Furthermore, it is simple to specify the set $\mathcal{H}(k)$ in order to enforce traditional performance specifications for regulation and tracking problems over finite horizons. For example, let system \eqref{eq:systemeq1}-\eqref{eq:systemeq2} be the discrete-time description of a continuous-time system over a horizon $[0,T]$. Then, suppose it is desired that the i-th output variable, $y_i$, initiating at $y_i(0)$, reaches a set-point, $y_i^{sp}$, within a settling-time, $t_s$, and with a steady-state error of $\lambda_s$. Moreover, suppose this should occur within a peak over-shoot of $y_i^p$ and that the output variable rises to within $\lambda_r$ of $y_i^{sp}$ with a rise-time of $t_r$. Letting the sampling time be denoted by $T_s$, one could use sampled versions of functions $\overline{h}_i(t)$ and $\underline{h}_i(t)$, like in Figure \ref{fig:perf_func}, where $t\in[0,T]$, and specify $Q_{j_1}(k) \equiv  C_i$, $\phi_{j_1}(k) = \overline{h}_i(kT_s)  + y_i^{sp}$, $Q_{j_2}(k) \equiv  -C_i$,  and $\phi_{j_2}(k) = -\underline{h}_i(kT_s)  - y_i^{sp}$, where $j_1$ and $j_2$ refer to two arbitrary rows of $Q(k)$.

\subsection*{Problem Statement}

Given the system \eqref{eq:systemeq1}-\eqref{eq:systemeq2} along with a time horizon, $K\geq 0$, a target set $\mathcal{X}_T$ and a time-varying target-tube, $\mathcal{H}(k)$, satisfying $\mathcal{X}_T=\mathcal{H}(K)$, find a linear time-varying feedback control, $u(k) = F(k)y(k)$, and an initial set, $\mathcal{X}_0\subset\mathcal{H}(0)$, such that the solution to \eqref{eq:systemeq1} satisfies $x(k)\in \mathcal{H}(k)$, for all $k \in  \{0,\dots, K\}$, for all $x_0\in\mathcal{X}_0$.

\section{Reachable Sets}\label{sec:reach_sets}

Reachable sets play a central role in the derivation of our results. In the next definition assume that a feedback matrix, $F(k)$, has been specified, resulting in the control $u(k) = F(k)y(k) = F(k)Cx(k)$.
\begin{definition}[One-step reachable set]
	Consider the closed-loop system:
	\begin{equation}
	x(k+1) = (\bm{A}(k) + \bm{B}(k)F(k)C )x(k),\label{eq:basic}
	\end{equation}
	for $k\in\mathbb{Z}$, with $[\bm{A}(k)\,\,\bm{B}(k)] \in \co\{[A_i\,\,B_i]\}_{i=1}^s$ and $x(k)\in S_1\subset\mathbb{R}^n$. The \emph{one-step reachable set} from $S_1$ via \eqref{eq:basic} is given by:
	\begin{align}
	\mathcal{R}&([\bm{A}(k) + \bm{B}(k)F(k)C],S_1) \nonumber\\
	&\triangleq \{x\in\mathbb{R}^n:x = \left(\bm{A}(k) + \bm{B}(k)F(k)C\right)x(k),\nonumber\\
	& \quad \quad x(k)\in S_1, [\bm{A}(k)\,\,\bm{B}(k)] \in \co\{[A_i\,\,B_i]\}_{i=1}^s\}. \nonumber
	\end{align}
\end{definition}
\vspace{-5mm}
The one-step reachable set is the set of all possible successor states. Note that $\mathcal{R}([\bm{A}(k) + \bm{B}(k)F(k)C],S_1) = \co\{\mathcal{R}([A_i + B_iF(k)C],S_1)\}.$

\begin{figure}[thpb]
	\begin{center}
		\includegraphics[width=0.85\columnwidth]{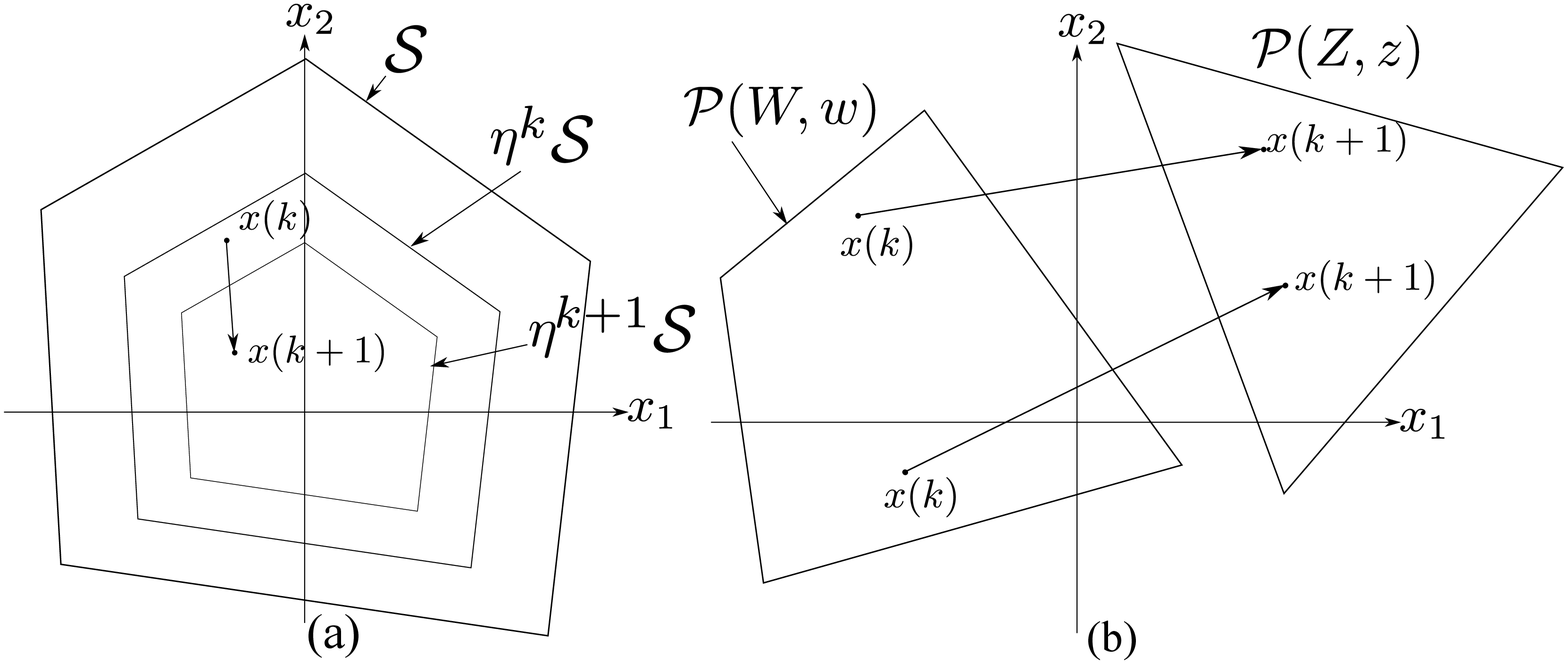}
		\caption{The figure labelled (a) summarises Proposition \ref{thm:original}: if its conditions are met, then at an arbitrary $k$ we have $x(k)\in \eta^k\mathcal{S}$, $0\leq\eta <1$. The figure labelled (b) summarises Theorem~\ref{thm:main}: if its conditions are met, then for $x(k)\in\mathcal{P}(W,w)$ it is guaranteed that $x(k+1)\in\mathcal{P}(Z,z)$, where $\mathcal{P}(W,w)$ and $\mathcal{P}(Z,z)$ are arbitrary polyhedral subsets of $\mathbb{R}^n$.}\label{fig:thm_adaptation}
	\end{center}
\end{figure}

In order to give the reader some background as to where our main result (Theorem~\ref{thm:main}) comes from, we briefly present ideas concerning polyhedral invariant sets. The next proposition, which we have slightly modified, originally appeared in \cite{Bitsoris:1988} and \cite{Bitsoris:Jun1988} in a form that guarantees invariance, i.e. the case where $\eta = 1$.
\begin{proposition}[ \cite{Bitsoris:1988}, \cite{Bitsoris:Jun1988}, \cite{Blanchini:2015}]\label{thm:original}
	Consider the linear system:
	\begin{equation}
	x(k+1) = Ax(k), \label{eq:no-cont-no-dist}
	\end{equation}
	$k\in\mathbb{Z}_{\geq 0}$, and let $\mathcal{P}(W,\textbf{1})$ be a polyhedral C-set. The set $\mathcal{P}(W,\textbf{1})$ is $\eta$-contractive for the system \eqref{eq:no-cont-no-dist} if and only if there exists a matrix $G$ that satisfies: 
	\begin{align}
	&G \geq 0, \nonumber\\
	&GW = WA, \nonumber\\
	&G\textbf{1}\leq \eta\textbf{1}, \nonumber
	\end{align}
	where $0\leq\eta < 1$.
\end{proposition}
We omit the proof in the interest of conciseness but present its ideas in Figure \ref{fig:thm_adaptation}. Next we state a result that was originally introduced in \cite{ESTERHUIZEN_2017}. It is, in turn, an adaptation of the result in Proposition~\ref{thm:original}. In our adaptation the state evolves from an arbitrary polyhedral set to another. We also mention that this theorem has been utilised in the study of finite-time stability with polyhedral domains, see \cite{Esterhuizen_2018}.

\begin{theorem}\label{thm:main}
	Consider the system $x(k+1) = Ax(k)$ along with two polyhedral sets, $\mathcal{P}(M_1,\mu)\subset\mathbb{R}^n$ and $\mathcal{P}(M_2,\nu)\subset\mathbb{R}^n$, with $x(k)\in\mathcal{P}(M_1,\mu)$ for an arbitrary $k\in\mathbb{Z}$. The following holds: $\mathcal{R}(A,\mathcal{P}(M_1,\mu))\subset\mathcal{P}(M_2,\nu)$ if and only if there exists a matrix $G$ satisfying:
	\begin{align}
	G\geq & 0, \label{thm1:cond1}\\
	GM_1= &M_2 A,\label{thm1:cond2}\\
	G\mu\leq& \nu.\label{thm1:cond3}
	\end{align}
\end{theorem}
\vspace{-6mm}
\begin{pf}
	(if)
	We have $x(k)\in \mathcal{P}(M_1,\mu)$ for an arbitrary $k\in\mathbb{Z}$, thus $M_1x(k) \leq \mu$. If there exists a $G$ satisfying \eqref{thm1:cond1}-\eqref{thm1:cond3} then $M_2Ax(k) = GM_1x(k) \leq G\mu \leq \nu$. Hence we have $Ax(k)\in\mathcal{P}(M_2,\nu)$ and thus $\mathcal{R}(A,\mathcal{P}(M_1,\mu)) \subset \mathcal{P}(M_2,\nu)$.
	
	(only if)
	We have $x(k) \in \mathcal{P}(M_1,\mu)$ and $Ax(k) \in \mathcal{P}(M_2,\nu)$ for an arbitrary $k\geq 0$. For each $j \in \{1,\dots,q_2 \}$, where $M_2\in\mathbb{R}^{q_2\times n}$, consider the following linear program:
	\begin{equation}
	\rho_j = \max_{x} M_2^{[j]}Ax,\,\,\, \mathrm{s.t.}\,\,\, M_1x \leq \mu.\nonumber
	\end{equation}
	The dual problem is
	\begin{equation}
	\rho_j = \min_{g_j} g_j^T \mu,\,\,\, \mathrm{s.t.}\,\,\, g_j^TM_1 = M_2^{[j]}A,\,\,\, g_j^T \geq 0,\nonumber
	\end{equation}
	where $g_j$ is the vector of dual variables for the $j-th$ linear program. Let $G^{[j]}$ (a row vector) be a feasible solution to the dual problem, and let $G$ be the matrix formed by stacking these $q_2$ solutions. Then, $GM_1 = M_2A$ and every entry of $G$ is nonnegative, so that \eqref{thm1:cond1} and \eqref{thm1:cond2} are true. Noting that $\rho_j = \max_{x\in\mathcal{P}(M_1,\mu)} M_2^{[j]} Ax \leq \nu_j$ and the fact that $\rho_j = G^{[j]}\mu$ (by strong duality), we have $G^{[j]}\mu \leq \nu_j$ for each $j$, which implies $G\mu\leq \nu$, which is \eqref{thm1:cond3}.
\end{pf}

We now generalise Theorem~\ref{thm:main} to the case of the uncertain system \eqref{eq:basic}.
\begin{corollary}\label{cor1}
	Consider the closed loop system \eqref{eq:basic} with an input $u(k) = F(k)y(k) = F(k)Cx(k)$, along with two polyhedral sets, $\mathcal{P}(M_1,\mu)\subset\mathbb{R}^n$ and $\mathcal{P}(M_2,\nu)\subset\mathbb{R}^n$, with $x(k)\in\mathcal{P}(M_1,\mu)$ for an arbitrary $k\in\mathbb{Z}$. The following holds: $\mathcal{R}([\bm{A}(k) + \bm{B}(k)F(k)C],\mathcal{P}(M_1,\mu))\subset\mathcal{P}(M_2,\nu)$ if and only if for every $i\in\{1,2,\dots,s\}$ there exists a matrix $G_i$ satisfying:
	\begin{align}
	G_i\geq & 0, \label{cor1:cond1}\\
	G_iM_1= &M_2 [A_i + B_iF(k)C],\label{cor1:cond2}\\
	G_i\mu\leq& \nu.\label{cor1:cond3}
	\end{align}
\end{corollary}
\vspace{-6mm}
\begin{pf}
	Because $\mathcal{P}(M_2,\nu)$ is convex, we have $\mathcal{R}([A_i + B_iF(k)C],\mathcal{P}(M_1,\mu))\subset\mathcal{P}(M_2,\nu)$ for $i\in\{1,2,\dots,s\}$ if and only if $\co\{\mathcal{R}([A_i + B_iF(k)C],\mathcal{P}(M_1,\mu))\}_{i=1}^s = \mathcal{R}([\bm{A}(k) + \bm{B}(k)F(k)C],\mathcal{P}(M_1,\mu)\subset\mathcal{P}(M_2,\nu)$. From Theorem~\ref{thm:main} we have $\mathcal{R}([A_i + B_iF(k)C],\mathcal{P}(M_1,\mu))\subset\mathcal{P}(M_2,\nu)$ for every $i\in\{1,2,\dots,s\}$ if and only if conditions \eqref{cor1:cond1}-\eqref{cor1:cond3} hold for $i\in\{1,2,\dots,s\}$, which completes the proof.
\end{pf}

\begin{remark}
	Corollary~\ref{cor1} says that for the system \eqref{eq:basic} one only needs to be concerned with the ``extreme'' matrices $[A_i\quad B_i]$.
\end{remark}

\section{Algorithm}\label{sec:algorithm}
We now present an algorithm that takes advantage of Corollary~\ref{cor1} to construct a time-varying feedback that addresses the problem statement. We introduce the set:
\[
\mathcal{X}(k) = \mathcal{P}(Q(k),\psi(k)),\,\,k\in\{0,1,\dots,K\},
\]
with $\mathcal{X}(K) = \mathcal{X}_T$, where $Q(k)$ is the same matrix that appears in the definition of $\mathcal{H}(k)$. Taking a dynamic programming approach, as in \cite{BERTSEKAS_1971}, the idea is to start with the target set, $\mathcal{X}_T$, and, working \emph{backwards} in time, to find a feedback $F(K-1)$ and a set $\mathcal{X}(K-1)\subset\mathcal{H}(K-1)$ such that $\mathcal{R}([\bm{A}(K-1) + \bm{B}(K-1)F(K-1)C],\mathcal{X}(K-1))\subset\mathcal{H}(K)$. We then find an $F(K-2)$ and an $\mathcal{X}(K-2)\subset\mathcal{H}(K-2)$ such that $\mathcal{R}([\bm{A}(K-2) + \bm{B}(K-2)F(K-2)C],\mathcal{X}(K-2))\subset\mathcal{X}(K-1)$, and continue going backwards in this fashion until $k=0$.
\vspace{-2mm}
\subsection*{Algorithm 1}
\begin{algorithmic}[1]
	\Statex \textbf{Inputs}: 
	$Q(k)$, $\phi(k)$ for $k = 0,\dots, K$
	\Statex \textbf{Outputs}: $\psi(k)$, $F(k)$ for $k = 0,\dots K-1$.
	\Statex \bf{Begin.}
	\Statex Let $\psi(K) \leftarrow \phi(K)$.
	\For{$k = K-1,K-2,\dots,0$}
	\State Solve LP1.
	\State Let $(\bar{\ee}(k),\{\bar{G}_i(k)\}_{i=1}^s,\bar{F}(k))$ be the solution to LP1.
	\State Let $F(k) \leftarrow \bar{F}(k)$.
	\If{$\bar{\ee}(k) = 0$}
	\State Let $\psi(k) \leftarrow \phi(k)$.
	\Else
	\State Solve LP2, with $\{\bar{G}_i(k)\}_{i=1}^s$.
	\State Let $\bar{\psi}(k)$ be the solution to LP2.
	\State Let $\psi(k) \leftarrow \bar{\psi}(k)$.
	\EndIf
	\EndFor
	\Statex \bf{End.}
\end{algorithmic}\label{alg:1}
\begin{LP}(LP1)
	\begin{align}
	&\min_{\ee(k),\{G_i(k)\}_{i=1}^s,F(k)}\,\,\,\,\, \ee(k)^T\bm{1}\nonumber\\
	&\mathrm{s.t.}\quad 
	\emph{\textrm{for }} i\in\{1,2,\dots,s\}:\nonumber\\
	&G_i(k) \geq 0, \nonumber\\
	&G_i(k)Q(k)
	= Q(k+1)[A_i + B_iF(k)C], \label{op1_const_2}\\
	& G_i(k)\phi(k)\leq \psi(k+1) + \ee(k), \label{op1_const_3}\\
	& \ee(k)\geq 0.\,\,\nonumber
	\end{align}
\end{LP}
\begin{LP}(LP2)
	\begin{align}
	\max_{\psi(k)}\,\,\,\,\, &\psi(k)^T\bm{1}\nonumber\\
	&\mathrm{s.t.}\quad 
	 \emph{\textrm{for }} i\in\{1,2,\dots,s\}:\nonumber\\
	& \bar{G}_i(k)\psi(k)\leq \psi(k+1), \label{op2_const_3}\\
	& \psi(k)\leq \phi(k).\label{op2_const_5}
	\end{align}
\end{LP}

\subsection{Explanation of Algorithm}
For an arbitrary $k\in\{0,1,\dots,K-1\}$ LP1 considers the set $\mathcal{X}(k+1)$, which satisfies $\mathcal{X}(k+1)\subset\mathcal{H}(k+1)$, along with the set $\mathcal{H}(k)$, and attempts to find an $F(k)$ such that:
\[
\mathcal{R}([\bm{A}(k) + \bm{B}(k)F(k)C(k)],\mathcal{H}(k))\subset\mathcal{X}(k+1).
\]
If $\bar{\ee}(k) \neq 0$ one can conclude that there does \emph{not} exist an $F(k)$ such that this is possible. However, there might still be an $\mathcal{X}(k)\subset\mathcal{H}(k)$ such that $\bar{F}(k)$, computed from LP1, results in:
\[
\mathcal{R}([\bm{A}(k) + \bm{B}(k)\bar{F}(k)C(k)],\mathcal{X}(k))\subset\mathcal{X}(k+1).
\]
LP2 then attempts to find an $\mathcal{X}(k)$ that satisfies this. The constraint \eqref{op2_const_5} ensures that $\mathcal{X}(k)\subset\mathcal{H}(k)$, and the cost in LP2 is included in order to find a ``large'' set contained in $\mathcal{H}(k)$. Thus, the idea is that at an arbitrary $k$ LP1 finds an $F(k)$ such that the resulting one-step reachable set is as close as possible to being contained in $\mathcal{H}(k+1)$. Then, LP2 tries to find a subset of $\mathcal{H}(k)$ such that this is true.

\begin{remark}
	The main difficulty associated with the conditions of Corollary~\ref{cor1} is that, if one assumes $G_i(k)$ and $\phi(k)$ are decision variables, the constraint \eqref{op1_const_3} is nonlinear and nonconvex. The algorithm is an attempt at working around this difficulty.
\end{remark}

\subsection{Solution to the Problem}

\begin{proposition}
	Suppose that Algorithm 1 successfully executes and produces a pair $(\psi(k),F(k))$ for $k\in\{0,1,\dots,K-1\}$. Then, for every $x_0\in\mathcal{X}_0 \triangleq \mathcal{X}(Q(0),\psi(0))$ the computed feedback results in a solution to the system \eqref{eq:systemeq1}-\eqref{eq:systemeq2} satisfying $x(k)\in\mathcal{X}(k)\subset\mathcal{H}(k),\,\,k\in\{0,1,\dots,K\}.$
\end{proposition}
\vspace{-6mm}
\begin{pf}
	The proof is by induction. Consider an arbitrary $k\in\{0,1,\dots,K-1\}$, with $x(k)\in\mathcal{X}(k)\subset\mathcal{H}(k)$. If the pair $(\psi(k),F(k))$ was produced from the solution of LP1 and $\bar{\ee}(k) = 0$, we can conclude from Theorem~\ref{thm:main} that:
	\[
	\mathcal{R}([\bm{A}(k) + \bm{B}(k)F(k)C],\mathcal{H}(k))\subset\mathcal{X}(k+1)\subset\mathcal{H}(k+1).
	\]
	Otherwise, if the pair $(\psi(k),F(k))$ was produced from the solution of LP2 (due to the fact that LP1 produced a solution with $\bar{\ee}(k) \neq 0$), we can conclude from Theorem~\ref{thm:main} that:
	\[
	\mathcal{R}([\bm{A}(k) + \bm{B}(k)F(k)C],\mathcal{X}(k))\subset\mathcal{X}(k+1)\subset\mathcal{H}(k+1).
	\]
	Thus, for $x(k)\in\mathcal{X}(k)$ at an arbitrary $k\in\{0,1,\dots,K-1\}$, it is guaranteed that $x(k+1)\in\mathcal{X}(k+1)\subset\mathcal{H}(k+1)$. We have $x_0\in\mathcal{X}_0$, and so the proof follows by induction.
\end{pf}

\subsection{Existence of Solution}\label{subsec:existence}
\noindent We introduce the following two assumptions:
\begin{description}
	\item[(A2)] For every $k\in\{0,1,\dots,K\}$ the target tube $\mathcal{H}(k)$ contains the origin in its interior.
	\item[(A3)] For an arbitrary $k$ the constraint $\psi(k)\geq 0$ is added to LP2.
\end{description}

\begin{proposition}\label{prop2}
	Consider the system \eqref{eq:systemeq1} - \eqref{eq:systemeq2}. Suppose (A1), (A2) and (A3) hold. Then there always exists a solution $(\psi(k),F(k))$, for $k\in\{0,1,\dots, K-1\}$, as the output of  Algorithm~1.
\end{proposition}
\vspace{-5mm}
\begin{pf}
	Because the entries of the matrices $A_i$ and $B_i$, $i\in\{1,2\dots,s\}$ are real (and thus finite) and the set $\mathcal{H}(k)$ is bounded for $k\in\{0,1,\dots,K\}$ (from (A1)), the one-step reachable set, $\mathcal{R}([\bm{A}(k) + \bm{B}(k)F(k)C],\mathcal{H}(k))$ for $k\in\{0,1,\dots,K-1\}$ is bounded for any choice of feedback matrix $F(k)$. Thus, there always exists a bounded polyhedral set, with the normals of its constituent hyperplanes specified by $Q(k+1)$, that contains this reachable set, which, via Corollary~\ref{cor1}, implies that there must exist a solution to LP1. (The vector $\ee(k)$ in LP1 may be chosen arbitrarily large, but finite.) The problem specified by LP2, with Assumption (A3), always has a solution: because of assumption (A2), for any $k$ we have $\phi(k)\geq 0$. Also, because of lines 6 and 10 of the algorithm and (A3), we have $\psi(k+1)\geq 0$. Thus, for any $k$, $\psi(k) = 0$ is always a solution to LP2 with Assumption (A3).
\end{pf}

\begin{remark}
	Assumptions (A2) and (A3) imply that one may always resort to the singleton $\{0\}$ as the set $\mathcal{X}(k)$, because
	$\mathcal{R}([\bm{A}(k) + \bm{B}(k)F(k)C],\{0\}) = \{0\} \in \mathcal{H}(k+1)$ for any feedback $F(k)$. For some problems, especially ones with large uncertainty, it may be that the sets $\mathcal{X}(k)$ become small neighbourhoods of the origin as $k$ iterates down from $K-1$. In this case, it may make sense to drop (A2) and (A3) at the cost of not having a guarantee of finding a solution from Algorithm~1.
\end{remark}

\section{Generalisations}\label{sec:generalisations}

\subsection{Additive Disturbance Term}

Consider the system:
\begin{align}
x(k + 1)  & =   \bm{A}(k)x(k) + \bm{B}(k)u(k) + 
Dv(k),\label{eq_sys_dist_1}\\
y(k) & = Cx(k), \label{eq_sys_dist_2}
\end{align}
where for every $k$, $\bm{A}(k)$ and $\bm{B}(k)$ are the same as before, and the disturbance, $v(k)\in\mathbb{R}^p$, is assumed to be located in a time-varying polyhedral set, given by
\[
\mathcal{V}(k) = \mathcal{P}(W(k), \gamma(k))\quad k \in\{0,\dots,K-1\},
\]
with $W(k)\in\mathbb{R}^{q_v(k)\times p}$, $q_v(k)\in \mathbb{Z}_{\geq 0}$, and $\gamma(k)$ a vector of compatible dimension for every $k$. At an arbitrary $k\in\mathbb{Z}$, consider the extended vector $(x(k)^T,v(k)^T)^T\in\mathbb{R}^{n + p}$. With a feedback, $u(k) = F(k)y(k)$, we can consider the equation:
\[
x(k+1) = [\bm{A}(k) + \bm{B}(k)F(k)C\quad D][x(k)^T,v(k)^T]^T,
\]
and define a modification of the one-step reachable set.
\begin{definition}
	Consider the closed-loop system:
	\begin{equation}
	x(k+1) = [\bm{A}(k) + \bm{B}(k)F(k)C]x(k) + Dv(k),\label{eq:basic_dist}
	\end{equation}
	for $k\in\mathbb{Z}$, with $[\bm{A}(k)\,\,\bm{B}(k)] \in \co\{[A_i\,\,B_i]\}_{i=1}^s$, $x(k)\in S_1\subset\mathbb{R}^n$ and $v(k)\in S_2\subset\mathbb{R}^p$. The \emph{one-step reachable set with disturbance} from $S_1$ and $S_2$ via \eqref{eq:basic_dist} is given by:
	\begin{align}
	\mathcal{R}&([\bm{A}(k) + \bm{B}(k)F(k)C\quad D],S_1, S_2) \nonumber\\
	&\triangleq \{x:x = [\bm{A}(k) + \bm{B}(k)F(k)C\quad D][x(k)^T,v(k)^T]^T,\nonumber\\
	& \quad x(k)\in S_1, v(k)\in S_2, [\bm{A}(k)\,\,\bm{B}(k)] \in \co\{[A_i\,\,B_i]\}_{i=1}^s\}. \nonumber
	\end{align}
\end{definition}
\noindent Note that:
\begin{align}
&\mathcal{R}([\bm{A}(k) + \bm{B}(k)F(k)C\quad D],S_1, S_2)\nonumber\\
& = \co\{\mathcal{R}([A_i + B_iF(k)C],S_1)\}\oplus DS_2\nonumber\\
& = \co\{\mathcal{R}([A_i + B_iF(k)C\quad D],S_1,S_2)\}.\nonumber
\end{align}
We generalise Corollary~\ref{cor1} as follows:
\begin{corollary}\label{cor2}
	Consider the system \eqref{eq:basic_dist}, along with three polyhedral sets, $\mathcal{P}(M_1,\mu)\subset\mathbb{R}^n$, $\mathcal{P}(M_2,\nu)\subset\mathbb{R}^n$,  and $\mathcal{P}(M_3,\xi)\subset\mathbb{R}^p$, with $x(k)\in\mathcal{P}(M_1,\mu)$ and $v(k)\in\mathcal{P}(M_3,\xi)$ for an arbitrary $k\in\mathbb{Z}$. The following holds: $\mathcal{R}([\bm{A}(k) + \bm{B}(k)F(k)C\quad D],\mathcal{P}(M_1,\mu),\mathcal{P}(M_3,\xi))\subset\mathcal{P}(M_2,\nu)$, if and only if for every $i\in\{1,2,\dots,s\}$ there exists a matrix $G_i$ satisfying:
	\begin{align}
	G_i\geq & 0, \label{cor2:eq1}\\
	G_i
	\left[
	\begin{array}{cc}
	M_1 & \bm{0}\\
	\bm{0} &M_3
	\end{array}
	\right] = &M_2 [A_i + B_iF(k)C\quad D],\label{cor2:eq2}\\
	G_i
	\left[
	\begin{array}{c}
	\mu\\
	\xi
	\end{array}
	\right]
	\leq& \nu.\label{cor2:eq3}
	\end{align}
\end{corollary}
\begin{pf}
	The proof of Theorem~\ref{thm:main} is easily adapted to the current setting, and we can conclude that $\mathcal{R}([A_i + B_iF(k)C\quad D],\mathcal{P}(M_1,\mu), \mathcal{P}(M_3,\xi))\subset \mathcal{P}(M_2,\nu)$ for every $i$ if and only if conditions \eqref{cor2:eq1}-\eqref{cor2:eq3} hold for every $i$. Then, as argued in the proof of Corollary~\ref{cor1}, $\mathcal{R}([A_i + B_iF(k)C\quad D],\mathcal{P}(M_1,\mu), \mathcal{P}(M_3,\xi))\subset \mathcal{P}(M_2,\nu)$ for every $i$ if and only if 
	\begin{align*}
	\mathcal{R}&([\bm{A}(k) + \bm{B}(k)F(k)C\quad D],\mathcal{P}(M_1,\mu), \mathcal{P}(M_3,\xi))\\
	&\subset \mathcal{P}(M_2,\nu),
	\end{align*}
	which completes the proof.
\end{pf}
\vspace{-5mm}
To now use Algorithm 1 to find a feedback sequence for the system with disturbance, the conditions \eqref{op1_const_2} and \eqref{op1_const_3} of LP1 need to replaced with:
\begin{equation}
G_i
\left[
\begin{array}{cc}
Q(k) & \bm{0}\\
\bm{0} & W(k)
\end{array}
\right] = Q(k+1) [A_i + B_iF(k)C\quad D],\label{eq:dist_mod_1}
\end{equation}
and
\begin{equation}
G_i(k)\left[\begin{array}{c}
\phi(k)\\
\gamma(k)
\end{array}\right]
\leq \psi(k+1) + \ee(k),\label{eq:dist_mod_2}
\end{equation}
respectively, and the condition \eqref{op2_const_3} of LP2 needs to be replaced with:
\begin{equation}
\bar{G}_i(k)\left[\begin{array}{c}
\psi(k)\\
\gamma(k)
\end{array}\right]\leq \psi(k+1).\label{eq:dist_mod_3}
\end{equation}

The existence result of Proposition~\ref{prop2} can be modified as follows; we need to enforce additional assumptions:
\begin{description}
	\item[(A4)] $D\mathcal{V}(k)\subset\mathcal{H}(k+1)$ for every $k\in\{0,1,\dots,K-1\}$,
	\item[(A5)] The constraint $D\mathcal{V}(k)\subset\mathcal{X}(k+1)$ is added to LP2 for $k\in\{K-2,K-3,\dots,0\}$.
\end{description}
Then, we can state:
\begin{proposition}\label{prop3}
	Consider the system \eqref{eq_sys_dist_1} - \eqref{eq_sys_dist_2}. Suppose (A1)-(A5) hold, that the conditions \eqref{op1_const_2} and \eqref{op1_const_3} of LP1 have been replaced with \eqref{eq:dist_mod_1} and \eqref{eq:dist_mod_2} respectively, and that constraint \eqref{op2_const_3} of LP2 has been replaced with \eqref{eq:dist_mod_3}. Then there always exists a solution $(\psi(k),F(k))$, for $k\in\{0,1,\dots, K-1\}$, as the output of Algorithm~1.
\end{proposition}
\begin{pf}
	As in the proof of Proposition~\ref{prop2}, because of Assumption (A1) there always exists a solution to LP1 with the conditions \eqref{op1_const_2} and \eqref{op1_const_3} replaced with conditions \eqref{eq:dist_mod_1} and \eqref{eq:dist_mod_2}, respectively. Now, let $k=K-1$ and suppose LP1 has been solved with $\ee(K-1)\neq 0$, producing $\bar{F}(K-1)$. Consider LP2 with the constraint \eqref{op2_const_3} replaced with \eqref{eq:dist_mod_3}, and $k=K-1$. Because of (A2)-(A4), with $\psi(K-1) = 0$, we have:
	\begin{align}
	&\mathcal{R}([\bm{A}(K-1) + \bm{B}(K-1)\bar{F}(K-1)C\quad D],\nonumber\\
	&\quad\quad\quad\quad\quad\quad\quad\quad\quad\quad\quad\quad\{0\}, \mathcal{P}(W(K-1),\gamma(K-1))\nonumber\\
	&= \{0\}\oplus D\mathcal{P}(W(K-1),\gamma(K-1))\subset \mathcal{H}(K).\nonumber
	\end{align}
	Otherwise, with $k\in\{K-2,K-3,\dots,0\}$ and Assumptions (A4) and (A5), we have with $\psi(k) = 0$:
	\begin{align}
	&\mathcal{R}([\bm{A}(k) + \bm{B}(k)\bar{F}(k)C\quad D],\{0\}, \mathcal{P}(W(k),\gamma(k))\nonumber\\
	&= \{0\}\oplus D\mathcal{P}(W(k),\gamma(k))\subset \mathcal{X}(k+1)\subset\mathcal{H}(k+1).\nonumber
	\end{align}
	Thus, under Assumptions (A2)-(A5), $\psi(k) = 0$ is always a solution to LP2 with \eqref{op2_const_3} replaced with \eqref{eq:dist_mod_3}.
\end{pf}

\begin{remark}
	If, for every $k$, the set $\mathcal{V}(k)$ is bounded and its vertices, denoted by $\mathsf{vert}(\mathcal{V}(k))$, are known, then imposing the constraints in (A5) on LP2 is straightforward. It translates to: $Q(k)Dv \leq \psi(k)$ for every $v\in\mathsf{vert}(\mathcal{V}(k))$.
\end{remark}

\subsection{Control Constraints}\label{subsec:control_constraints}

Considering the system \eqref{eq:systemeq1}-\eqref{eq:systemeq2}, or \eqref{eq_sys_dist_1}-\eqref{eq_sys_dist_2}, suppose the control is required to be constrained in a polyhedral set for every $k$, i.e.:
\begin{equation}
u(k) \in \mathcal{P}(U(k),\theta(k)).\label{constraint_control}
\end{equation}
If we introduce the assumption:
\begin{description}
	\item[(A6)] The vertices of the set $\mathcal{H}(k)$, denoted $\mathsf{vert}(\mathcal{H}(k))$, are known for every $k$,
\end{description}
then it is straightforward to show that the constraint \eqref{constraint_control} can be enforced by including the following constraint:
\[
U(k)F(k)Ch \leq \theta(k)\text{ for every }h\in\mathsf{vert}(\mathcal{H}(k)),
\]
in LP1. However, note that with this constraint it is not guaranteed that a solution will exist to LP1.

\section{Stability}\label{sec:stab}

Note that for system \eqref{eq_sys_dist_1}-\eqref{eq_sys_dist_2} it is impossible to achieve asymptotic stability because of the presence of a consistent disturbance. However, it may be desirable that the state remain in the target set for all $k\geq K$. Briefly, one way this may be achieved is as follows: consider the system \eqref{eq_sys_dist_1}-\eqref{eq_sys_dist_2} in closed loop with a particular (static) linear control law $u(k) = \hat{F}y(k)$, resulting in the closed-loop system $x(k+1) =  A_{cl}(k)x(k) + Dv(k),$
where $A_{cl}(k) = [\bm{A}(k) + \bm{B}(k)\hat{F}C]$. Assume that the disturbance set is time-invariant, i.e., $\mathcal{V}(k)\equiv \mathcal{V}$, and define the following set:
\begin{definition}[Robust Invariant set]
	A set $S\subset\mathbb{R}^n$ is said to be a \emph{robust invariant set} with respect to the closed loop system $x(k+1) =A_{cl}(k)x(k) + Dv(k)$ provided that for all $k$, for all $x\in S$ and for all $v\in\mathcal{V}$ we have $A_{cl}(k)x + Dv \in S$.
\end{definition}
The state can then be made to remain in $\mathcal{X}_T$ for all $k\geq K$ as follows: suppose the control law $u(k) = \hat{F}y(k)$ renders a set $S\subset\mathbb{R}^n$ robustly invariant with respect to the closed-loop system. Now, execute Algorithm~1 under the assumption:
\begin{description}
	\item[(A7)] $\mathcal{X}_T\subset S$,
\end{description}
to obtain the feedback sequence $\bar{F}(k)$. Then, clearly, the control law:
\[
u(k)=\begin{cases}
\bar{F}(k)y(k)\quad k\in\{0,1,\dots, K-1\},\\
\hat{F}y(k)\quad k\geq K,
\end{cases}
\]
results in $x(k)\in\mathcal{X}(k)$ for $k\in\{0,1,\dots,K\}$, and $x(k)\in\mathcal{X}_T$ for all $k\geq K$. This is an idea that shares some resemblance to ``dual mode'' control that appears in MPC stability studies, see for example \cite{Mayne2000}. We note that the effective computation of robustly invariant sets is a difficult problem in general, and the reader is directed to the works \cite{Blanchini:2015,Kerrigan2000RobustCS,Blanchini1990}, and the references therein.


\section{Example}\label{sec:example}

\begin{figure}[thpb]
	\begin{center}
		\includegraphics[width=0.8\columnwidth]{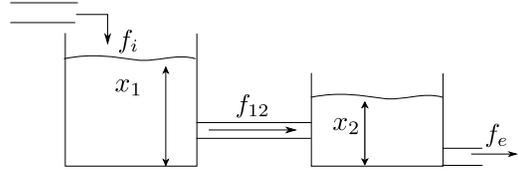}
		\caption{Two coupled water tanks. The control consists of $f_i$, the input flow rate, and $f_e$, the exit flow rate.}
		\label{fig8:tanks}
	\end{center}
\end{figure}

Consider the well-known example of two coupled tanks, as in Figure~\ref{fig8:tanks}. The left and right tanks, labelled tank 1 and tank 2, respectively, have cross-sectional areas of $R_1$~$\text{m}^2$ and $R_2$~$\text{m}^2$. Water can be pumped into tank 1 at a rate of $f_i$~$\text{m}^3.\text{s}^{-1}$ and out of tank 2 at a rate of $f_e$~$\text{m}^3.\text{s}^{-1}$. The two tanks are connected and water may freely flow between them, this flow being denoted by $f_{12}$.

\begin{figure}
	\centering
	\begin{tabular}{@{}c@{}}
		$\vspace{-3mm}$
		\includegraphics[width=0.9\columnwidth]{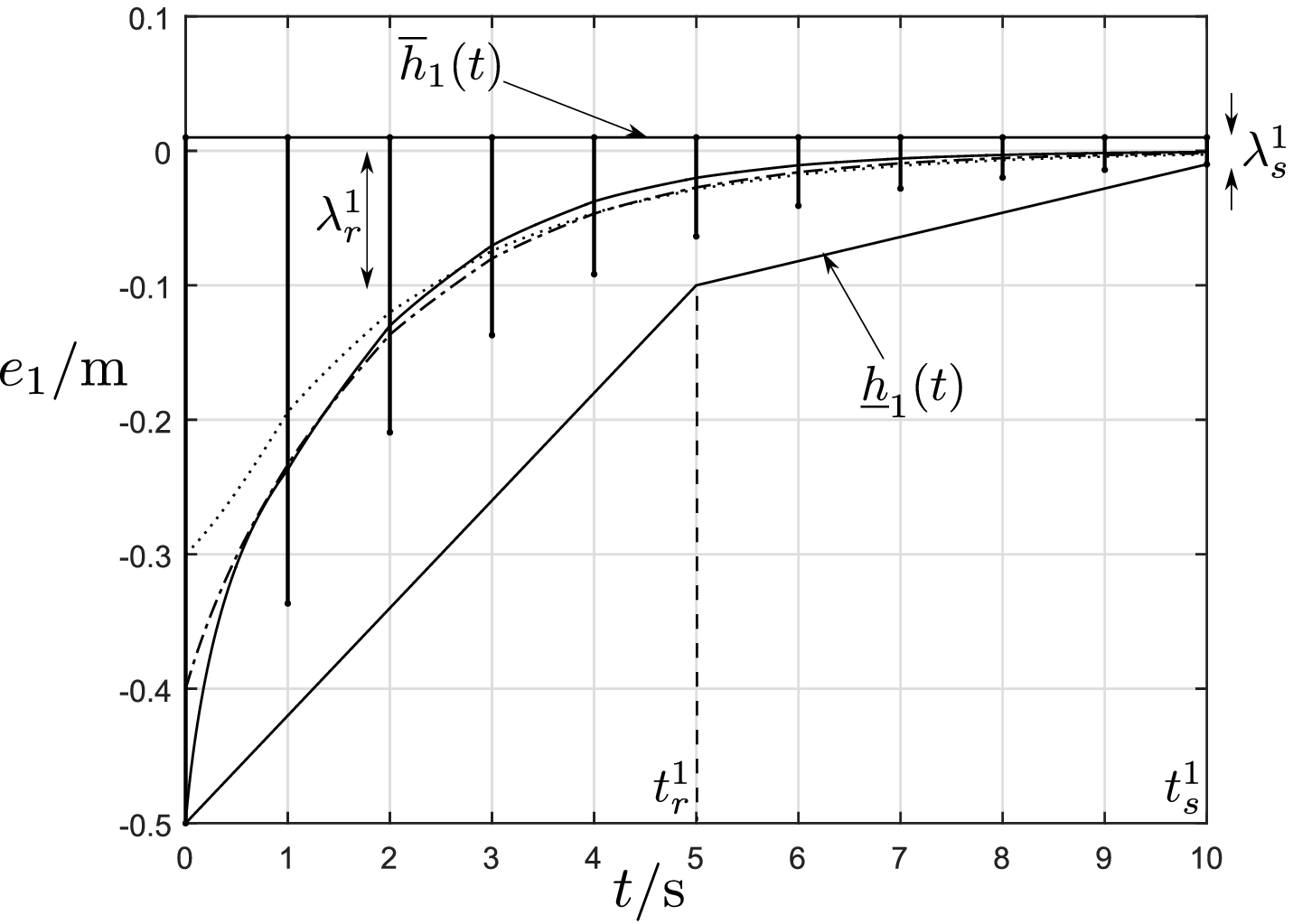}
		\small (a)
	\end{tabular}
	
	\vspace{\floatsep}
	
	\begin{tabular}{@{}c@{}}
		$\vspace{-3mm}$
		\includegraphics[width=0.9\columnwidth]{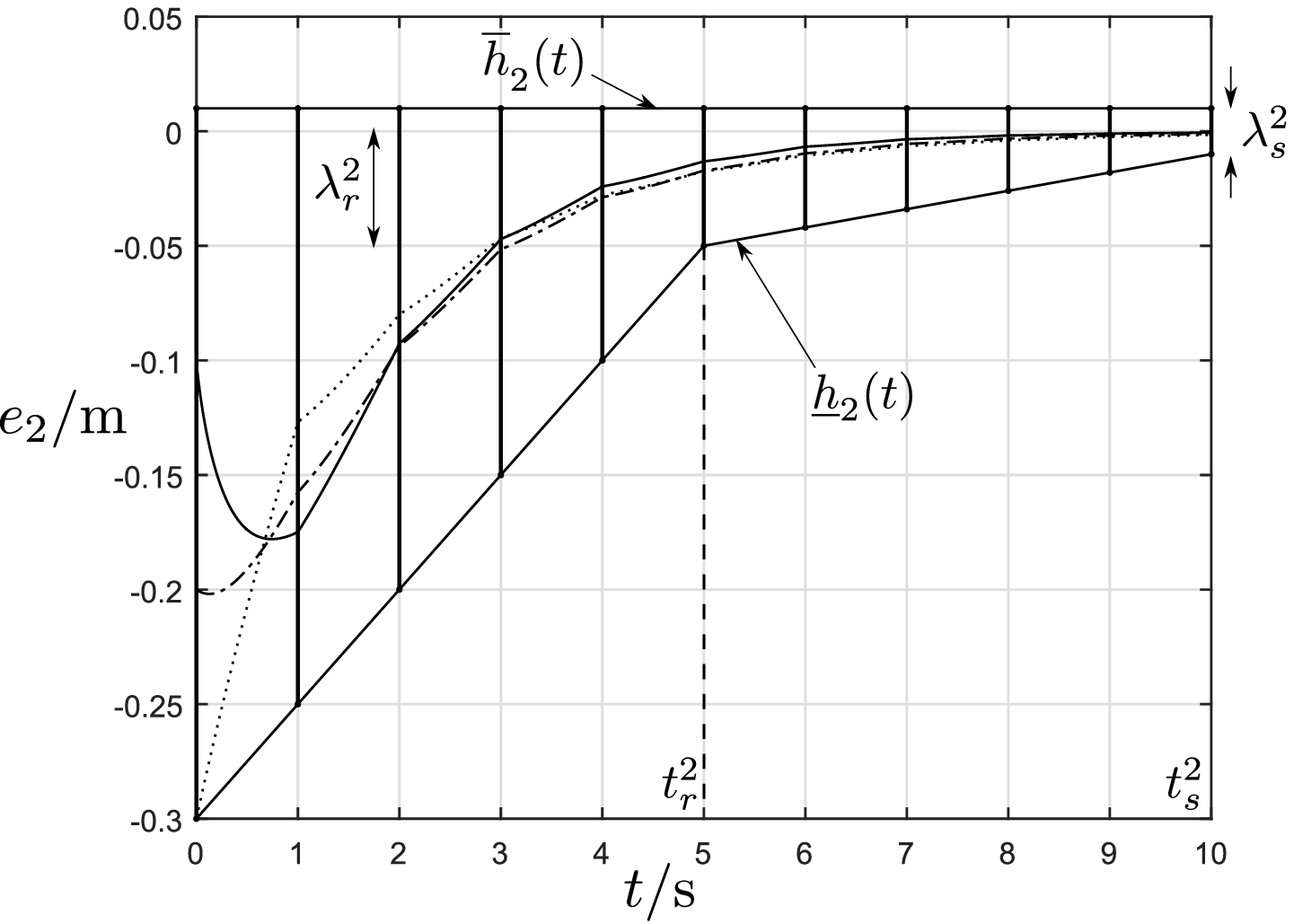}
		\small (b)
		$\vspace{3mm}$
	\end{tabular}
	
	\begin{tabular}{@{}c@{}}
		$\vspace{-3mm}$
		\includegraphics[width=0.9\columnwidth]{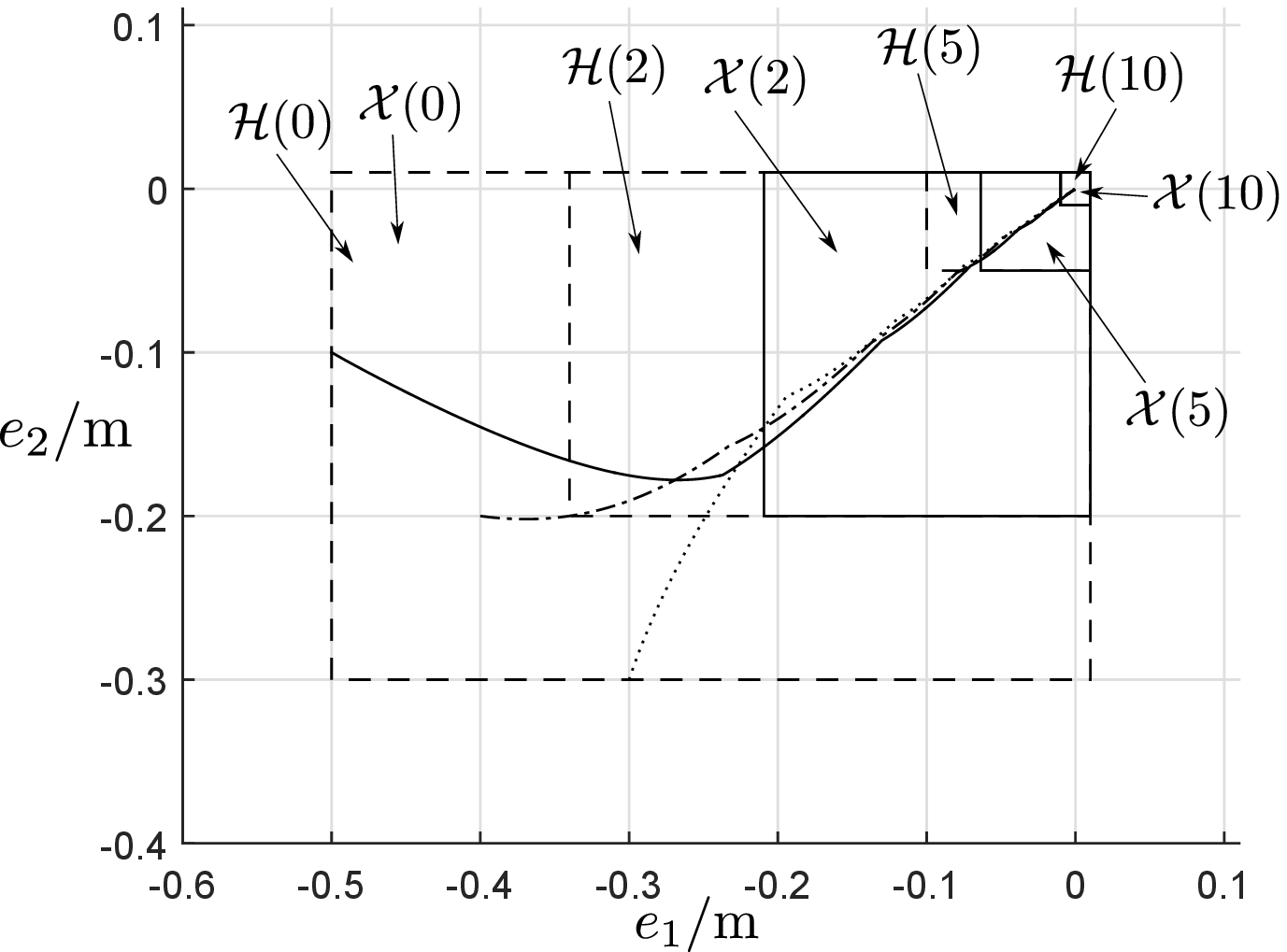}
		\small (c)
	\end{tabular}
	
	\caption{Results from the example in Section~\ref{sec:example}. The figure labelled (c) shows the target tube, $\mathcal{H}(k)$, with dashed boundaries, along with the computed set $\mathcal{X}(k)$, for some $k$'s. The figures labelled (a) and (b) show the projection of the set $\mathcal{X}(k)$ onto $e_1$ and $e_2$, respectively (the solid black vertical line segments). Three simulations of the \emph{nonlinear} system, each with a different value for $R_1$ and a different initial condition, $e_{\text{ini}}\in\mathcal{X}(0)$, are also shown ($R_1 = 3$ solid, $R_1 = 4$ dash-dotted, $R_1 = 5$ finely dashed.)}\label{fig:example}
\end{figure}

For completeness, we briefly show the derivation of a suitable linear model for this system, as done in \cite[Ch. 18]{Goodwin:2000}. Let $x_1$ and $x_2$ be the water level, in metres, in tank 1 and 2, respectively. If we define $u_1 \triangleq f_i/R_1$ and $u_2 \triangleq f_e/R_2$, a nonlinear model for this system is given by:
\begin{align}
&\dot{x}_1 = -L_1\sqrt{x_1 - x_2} + u_1, \nonumber\\
&\dot{x}_2 = L_2\sqrt{x_1 - x_2} + u_2, \nonumber
\end{align}
where $L_i = \frac{\sqrt{2g}}{R_i}$, $i=1,2$, $g = 10$ $\text{m}.\text{s}^{-2}$ is the acceleration due to gravity. Let $(\bar{x}_1,\bar{x}_2)$, with $\bar{x}_1 \geq \bar{x}_2$, be a desired state.  We introduce the error: $e_1(t) \triangleq x_1(t) - \bar{x}_1$, $e_2(t) \triangleq x_2(t) - \bar{x}_2$, and a new control: $\tilde{u}_1 \triangleq u_1 - L_1\sqrt{\bar{x}_1 - \bar{x}_2}$, $\tilde{u}_2 \triangleq u_2 + L_2\sqrt{\bar{x}_1 - \bar{x}_2}$. If we linearise about the equilibrium point $(e_1,e_2,\tilde{u}_1,\tilde{u}_2) = (0,0,0,0)$, we get the system:
\begin{equation}
\dot{e}(t) =
\frac{1}{2}\left( \bar{x}_1 - \bar{x}_2 \right)^{-\frac{1}{2}}
\left[
\begin{array}{cc}
-L_1& L_1\\
L_2 & -L_2
\end{array}
\right]e(t) 
+
\left[
\begin{array}{cc}
1 & 0\\
0 & 1
\end{array}
\right]\tilde{u}(t). \nonumber
\end{equation}
Our particular problem is as follows: both tanks have a height of 3 m, $R_2 = 5$, and $R_1\in\{3,4,5\}$. Our desired state is $(\bar{x}_1,\bar{x}_2) = (2, 1.6)$, and we want to find a time-varying linear feedback, along with a set of initial conditions, such that the following performance specifications are met for all initial conditions: $t^1_r = 5s$, $\lambda^1_r = 0.1$, $t^1_s = 10s$, $\lambda^1_s = 0.01$, $e_1$-overshoot = $0.01$; $t^2_r = 5s$, $\lambda^2_r = 0.05$, $t^2_s = 10s$, $\lambda^2_s = 0.01$, $e_2$-overshoot = $0.01$, where the superscripts 1 and 2 refer to tank 1 and 2 respectively. See Figure~\ref{fig:example} for further clarification. Moreover, we want to achieve this by \emph{only} measuring the water level in tank 2, i.e. $\tilde{u}(k) = F(k)e_2(k)$. We have three models, one for each value of $R_1$, which we discretise using a sampling time of $T_s = 1$ s. From our performance requirements we identify the target set $\mathcal{X}_T = [-0.01,0.01]\times [-0.01,0.01]$. The set $\mathcal{H}(k)$ is easily specified by referring to $\overline{h}_i(t)$ and $\underline{h}_i(t)$, $i=1,2$, in Figure~\ref{fig:example} (a) and (b). We specify the matrix $U(k)$ and $\theta(k)$, as explained in Subsection~\ref{subsec:control_constraints}, to ensure $f_i\geq 0$ and $f_e\leq 0$. We run Algorithm~1 to obtain the desired feedback along with the sets $\mathcal{X}(k)$, shown in Figure~\ref{fig:example}.

\section{Conclusion}\label{sec:conclusion}

In this paper we have presented a new ``target-tube'' approach to the design of a control law that guarantees desired transient performance for linear polytopic difference inclusions. A notable aspect of our study is that we produce a time-varying linear output-feedback law under the assumption of a polyhedral target-tube, $\mathcal{H}(k)$, which, to our knowledge, has not appeared in the literature. We presented a theorem that provides necessary and sufficient conditions for the state of the system to evolve from one polyhedral set to another. We then specified an algorithm that uses this theorem to construct a time-varying feedback along with a set, $\mathcal{X}(k)$, that satisfies $\mathcal{X}(k)\subset\mathcal{H}(k)$, such that the trajectory of the system, initiating at any point in $\mathcal{X}(0)$ is contained in $\mathcal{X}(k)$ for all $k$. We provided generalisations that allows one to include a bounded additive disturbance term and constraints on the input. We stated some assumptions under which it is guaranteed that the algorithm would successfully execute with a solution, noting that this may result in conservatism.

\begin{ack}                               
Qing-Guo WANG acknowledges the financial support of the National Research Foundation of South Africa (Grant Number: 113340) and Oppenheimer Memorial Trust grant, which partially funded his research on this work 
\end{ack}

\bibliographystyle{plain}        
\bibliography{autosam}           



\end{document}